\documentclass[11pt,reqno]{amsart}
\usepackage{color}
\usepackage{amssymb,latexsym}
\usepackage{verbatim}
\usepackage{graphicx}
\usepackage{url}		
\calclayout


\newcommand{\minus}{\scalebox{0.75}[1.0]{$-$}}

\theoremstyle{plain}
\numberwithin{equation}{section}
\newtheorem{thm}{Theorem}[section]
\newtheorem{theorem}[thm]{Theorem}

\newtheorem{definition}[thm]{Definition}

\begin{document}
\setcounter{page}{1}

\title{Fibonacci vector and matrix $p$-norms}
\author{Francisco Salas-Molina}
\address{Universitat Polit\`ecnica de Val\`encia, Alcoy, 03801, Spain \\}
\email{frasamo@upv.es}

\begin{abstract}
This paper delves into vector and matrix norms of Fibonacci numbers. Two classes of Fibonacci vectors and a parametric $p$-norm are defined. From this definition, several properties of Fibonacci vector and matrix $p$-norms are described by varying parameter $p$. A closed-form expression is given to obtain the value of $p$, setting the difference between the $p$-norm and the infinite norm below a given threshold. A new class of symmetric $k$-Fibonacci matrix is defined such that a simple reorganization simplifies the computation of its $p$-norm. The analysis is extended to $p$-distances when considering the norm of the difference of two vectors (matrices) of the same size. 
\end{abstract}

\maketitle

\section{Introduction\label{sec:intro}}
\belowdisplayskip=12pt
\abovedisplayskip=12pt

The Fibonacci sequence ${F_n}$ is defined as:
\begin{equation}
    F_1 = F_2 = 1 \; \text{and, for} \; n > 2, \; F_n = F_{n-1} + F_{n-2}.    
\end{equation}

$F_n$ is called the $n$-th Fibonacci number, and the Fibonacci sequence is:
\begin{equation}
    (F_0 := 0), 1, 1, 2, 3, 5, 8, 13, 21, 34, 55, 89, 144, \ldots
\end{equation}

A closed-form expression for the $n$-th Fibonacci number is given by the Binet formula \cite{koshy2019} involving the golden ratio $\varphi$:
\begin{equation}
    F_n =  \frac{\varphi^n + (-\varphi)^{-n}}{\sqrt{5}} = \frac{\varphi^n + (1-\varphi)^{n}}{\sqrt{5}}.
\end{equation}

The sum of the first $n$ Fibonacci numbers is given by \cite{benjamin2003}:
\begin{equation}
    \sum_{i=1}^n F_n = F_{n+2}-1.
\end{equation}

The sum of the squares of the first $n$ Fibonacci numbers is given by \cite{benjamin2003}:
\begin{equation}
    \sum_{i=1}^n F_n^2 = F_{n} \cdot F_{n+1}.
\end{equation}

Finally, the sum of the cubes of the first $n$ Fibonacci numbers is given by \cite{benjamin2009}:
\begin{equation}
    \sum_{i=1}^n F_n^3 = \frac{1}{2} (F_n \cdot F_{n+1}^2+(-1)^n \cdot F_{n-1} + 1).
\end{equation}

Generalized $k$-Fibonacci numbers were studied in \cite{falcon2007,falcon2007b,grossman1997,lee2003,lee2003b} with alternative definitions. The $k$-Fibonacci sequence $\{g(k)_n\}$ is defined in \cite{lee2003} as:
\begin{equation}
    g(k)_1 = \ldots = g(k)_{k-2} =  0, \; g(k)_{k-1} = g(k)_{k} = 1
\end{equation}
and for $n>k\geq 2$:
\begin{equation}
    g(k)_n = g(k)_{n-1} + g(k)_{n-2} + \ldots + g(k)_{n-k}
\end{equation}
that reduces to the Fibonacci sequence when $k=2$. From this $k$-Fibonacci sequence, two types of Fibonacci matrices were also introduced in \cite{lee2003}. First, an $n \times n$ $k$-Fibonacci matrix $\mathcal{F}(k)_n=\left[f(k)_{ij} \right]$ defined as:
\begin{equation}
    f(k)_{ij}= \left\{\begin{array}{ll}g_{i-j+1} & \; i-j+1 \geq 0 \\ 0 & \; i-j+1 < 0 . \end{array}\right.
\end{equation}

This definition places consecutive $k$-Fibonacci numbers in a lower triangular matrix form. An example of a $2$-Fibonacci matrix of order $5$ is:
\begin{equation}
\mathcal{F}(2)_5 = \left[ \begin{array}{lllll} 1 & 0 & 0 & 0 & 0 \\ 1 & 1 & 0 & 0 & 0 \\ 2 & 1 & 1 & 0 & 0 \\   3 & 2 & 1 & 1 & 0 \\  5 & 3 & 2 & 1 & 1 \end{array} \right].
\end{equation}

Second, an $n \times n$ symmetric $k$-Fibonacci matrix $\mathcal{Q}(k)_n=\left[q(k)_{ij} \right]$ defined as:
\begin{equation}
    q(k)_{ij}= q(k)_{ji} = \left\{\begin{array}{ll}  \sum_{l=1}^i g_i^2 & \; i=j \\ \sum_{l=1}^k q_{i,i-l} & \; i+1 \leq j \end{array}\right.
\end{equation}
where $g(k)_{ij}=0$ for $j \leq 0$. An example of a $2$-Fibonacci symmetric matrix of order $5$ is:
\begin{equation}
\mathcal{Q}(2)_5 = \left[ \begin{array}{rrrrr} 1 & 1 & 2 & 3 & 5 \\ 1 & 2 & 3 & 5 & 8 \\ 2 & 3 & 6 & 9 & 15 \\   3 & 5 & 9 & 15 & 24 \\  5 & 8 & 15 & 24 & 40 \end{array} \right].
\end{equation}

\section{Fibonacci vector $p$-norms}

Before discussing Fibonacci vector $p$-norms, let us first define norms:

\begin{definition}
    \textbf{Norm} \cite{horn2012}. Let $\boldsymbol{V}$ be an $n$-dimensional vector space over the field of real numbers $\mathbf{R}$. Function $||\cdot||: \boldsymbol{V} \rightarrow  \mathbf{R}$ is a norm if for all $\boldsymbol{x},\boldsymbol{y} \in \boldsymbol{V}$, and all $c \in \mathbf{R}$:
    \begin{itemize}
        \setlength\itemsep{1em}
        \item Non-negativity: $||\boldsymbol{x}|| \geq 0$.
        \item Positivity: $||\boldsymbol{x}|| = 0 \iff \boldsymbol{x} =  \boldsymbol{0}_n$, whwere $\boldsymbol{0}_n$ is an $n$-dimensional vector of zeros.
        \item Homogeneity: $||c \boldsymbol{x}|| = |c| \cdot ||\boldsymbol{x}||$.
        \item Triangle inequality: $||\boldsymbol{x}+\boldsymbol{y}|| \leq ||\boldsymbol{x}|| + ||\boldsymbol{y}||$.
    \end{itemize}
    \label{def:norm}
\end{definition}

A special class of norms is the $p$-norm, denoted by $||\boldsymbol{x}||_p$ for any vector $\boldsymbol{x}=(x_1, x_2, \ldots, x_n)$, and defined as follows: 
\begin{equation}
    ||\boldsymbol{x}||_p = \left( \sum_{i=1}^n |x_i|^p\right)^{1/p}.
\end{equation}

Consider now the set of an $n$-dimensional Fibonacci vector by extracting the first column vector of matrices $\mathcal{F}(2)_n$ or $\mathcal{Q}(2)_n$, denoted by $\boldsymbol{q}_n=(F_{1}, F_{2}, \ldots, F_{n})$. Then, its $p$-norm is:
\begin{equation}
    ||\boldsymbol{q}_n||_p= \left( \sum_{i=1}^n |F_i|^p\right)^{1/p}.
\end{equation}

In what follows, we derive several properties from basic identities described in Section \ref{sec:intro}:

\begin{enumerate}
    \setlength\itemsep{1em}
    \item $||\boldsymbol{q}_n||_p >0$ for $n \geq 1$
    \item $||\boldsymbol{q}_n||_1= \sum_{i=1}^n F_i = F_{n+2}-1$
    \item $||\boldsymbol{q}_n||_2^2= \sum_{i=1}^n F_i^2 = F_{n} \cdot F_{n+1}$
    \item $||\boldsymbol{q}_n||_3^3= \sum_{i=1}^n F_i^3 =\frac{1}{2} (F_n \cdot F_{n+1}^2+(-1)^n \cdot F_{n-1} - 1)$
    \item $||\boldsymbol{q}_n||_{\infty}= \displaystyle{\lim_{p \to \infty} \left( \sum_{i=1}^n F_i^p\right)^{1/p}} = max \{F_1, F_2, \ldots, F_n\} = F_n$
    \item $||\boldsymbol{q}_n||_{-\infty} = \displaystyle{\lim_{p \to -\infty} \left( \sum_{i=1}^n F_i^p\right)^{1/p}} = min \{F_1, F_2, \ldots, F_n\} = F_1$
    \item $||\boldsymbol{q}_n||_{-1}= \left(\sum_{i=1}^n 1/F_i\right)^{-1}=\left( \sqrt{5} \sum_{i=1}^n  \frac{1}{\varphi^{i} - (1-\varphi)^i}\right)^{-1}$ 
    \item $||\boldsymbol{q}_n||_0^n= \prod_{i=1}^n F_i$. \\
\end{enumerate}


By focusing on the value of $p$ that lays the difference between the $p$-norm and the infinite norm below a given threshold, we derive the following result:

\begin{theorem}
    Given $\boldsymbol{q}_n$, there exists some $p < \infty$ such that $||\boldsymbol{q}_n||_p - ||\boldsymbol{q}_n||_{\infty} \leq  \varepsilon$ for some value $\varepsilon>0$. The threshold value of $p$ is given by:
    \begin{equation}
          p \geq \frac{\ln{n}}{\ln{(\varepsilon/F_n+1)}}.
    \end{equation}
    \label{thm:epsilon}
\end{theorem}
\begin{proof}
    Through H\"older's inequality \cite{holder1889,maligranda1998}  in $\mathbb{R}^n$, when $0<p<s$, a $p$-norm relates to an $s$-norm through the following inequality:
\begin{equation}
   ||\boldsymbol{q}_n||_p \leq  n^{1/p -1/s} ||\boldsymbol{q}_n||_s.
   \label{eq:holder}
\end{equation}

Setting $s=\infty$ and subtracting $||\boldsymbol{q}_n||_{\infty}$ from both sides of the inequality:
\begin{equation}
    ||\boldsymbol{q}_n||_p - ||\boldsymbol{q}_n||_{\infty} \leq  n^{1/p} ||\boldsymbol{q}_n||_{\infty} - ||\boldsymbol{q}_n||_{\infty} \leq  \varepsilon
\end{equation}
\begin{equation}
    n^{1/p} F_n - F_n \leq  \varepsilon.
\end{equation}

Taking natural logarithms, we find the threshold value for $p$:
\begin{equation}
    \ln{n^{1/p}} \leq \ln{(\varepsilon/F_n+1)}
\end{equation}
\begin{equation}
    p \geq \frac{\ln{n}}{\ln{(\varepsilon/F_n+1)}}.
    \label{eq:pepsilon}
\end{equation}
\end{proof}

From inequality \eqref{eq:holder}, two additional properties are derived:
\begin{enumerate}
    \setlength\itemsep{1em}
    \setcounter{enumi}{8}
    \item $||\boldsymbol{q}_n||_p \geq ||\boldsymbol{q}_n||_s \; \forall p < s$
    \item $||\boldsymbol{q}_n||_p$ is decreasing in $p$.\\
\end{enumerate}


Let us focus now on $r$-dimensional Fibonacci vectors to define an $(n,r)$-Fibonacci vector as:
\begin{equation}
    \boldsymbol{q}_{n,r}=(F_{n+1}, F_{n+2}, \ldots, F_{n+r}).
\end{equation}

Along the lines of properties (1) to (10), different values of $p$ lead to the following properties for $n,r \geq 1$:
\begin{enumerate}
    \setlength\itemsep{1em}
    \setcounter{enumi}{10}
    \item $||\boldsymbol{q}_{n,r}||_p >0$ 
    \item $||\boldsymbol{q}_{n,r}||_1= ||\boldsymbol{q}_{n+r}||_1-||\boldsymbol{q}_n||_1 = F_{n+r+2} - F_{n+2}$ 
    \item $||\boldsymbol{q}_{n,r}||_2^2 = ||\boldsymbol{q}_{n+r}||_2^2  - ||\boldsymbol{q}_n||_2^2 = F_{n+r} \cdot F_{n+r+1} - F_{n} \cdot F_{n+1}$
    \item $||\boldsymbol{q}_{n,r}||_3^3 = ||\boldsymbol{q}_{n+r}||_3^3  - ||\boldsymbol{q}_n||_3^3 = \frac{1}{2} (F_{n+r} \cdot F_{n+r+1}^2+(-1)^{n+r} \cdot F_{n+r-1} - 1) - \frac{1}{2} (F_n \cdot F_{n+1}^2+(-1)^n \cdot F_{n-1} - 1)$
    \item $||\boldsymbol{q}_{n,r}||_{-1}= \left( \sqrt{5} \sum_{i=n+1}^{n+r} \frac{1}{\varphi^{i} - (1-\varphi)^i}\right)^{-1}$
    \item $||\boldsymbol{q}_{n,r}||_{\infty}= F_{n+r}$
    \item $||\boldsymbol{q}_{n,r}||_{-\infty}= F_{n+1}$
    \item $||\boldsymbol{q}_{n,r}||_0^n= \prod_{i=n+1}^{n+r} F_i$
    \item $||\boldsymbol{q}_{n,r}||_p \geq ||\boldsymbol{q}_{n,r}||_s \; \forall p < s$
    \item $||\boldsymbol{q}_{n,r}||_p$ is decreasing in $p$.\\
\end{enumerate}

\section{Fibonacci matrix $p$-norm}

Definition \ref{def:norm} also applies to an arbitrary $m \times n$ matrix $A=\left[a_{ij}\right]$, and its $p$-norm, denoted by $||A||_p$, can be defined as follows: 
\begin{equation}
    ||A||_p = \left(\sum_{i=1}^m \sum_{j=1}^n |a_{ij}|^p\right)^{1/p}.
\end{equation}

Note that when $p=2$, $||A||_2$ is the Frobenius matrix norm. Consider now the $2$-Fibonacci matrix of order $n$ as defined in the introduction \cite{lee2003}:

\begin{equation}
\mathcal{F}(2)_n = \left[ \begin{array}{ccccc} F_1 & 0 & 0 & 0 & 0 \\ F_2 & F_1 & 0 & 0 & 0 \\ F_3 & F_2 & F_1 & \ldots & 0 \\   \vdots & \vdots & \vdots & \ddots & \vdots \\  F_n & F_{n-1} & F_{n-2} & \ldots & F_1 \end{array} \right].
\end{equation}

Note that $\mathcal{F}(2)_n$ is built from the concatenation of $n$ vectors $\boldsymbol{q}_i$ with $i\in\{1,\ldots,n\}$. Then $p$-norm of $\mathcal{F}(2)_n$ is given by:
\begin{equation}
    ||\mathcal{F}(2)_n||_p^p = \sum_{i=1}^n ||\boldsymbol{q}_i||_p^p.
\end{equation}

Different values of $p$ lead to the following properties:

\begin{enumerate}
    \setlength\itemsep{1em}
    \setcounter{enumi}{20}
    \item $||\mathcal{F}(2)_n||_1 = \sum_{i=1}^n ||\boldsymbol{q}_i||_1 = \sum_{i=1}^n (F_{i+2}-1) = ||\boldsymbol{q}_{2,n+2}||_1 - n $
    \item $||\mathcal{F}(2)_n||_2^2 = \sum_{i=1}^n ||\boldsymbol{q}_i||_2^2= \sum_{i=1}^n F_{i} \cdot F_{i+1}$
    \item $||\mathcal{F}(2)_n||_3^3 = \sum_{i=1}^n ||\boldsymbol{q}_i||_3^3= \frac{1}{2} \sum_{i=1}^n (F_i \cdot F_{i+1}^2+(-1)^i \cdot F_{i-1} - 1)$
    \item $||\mathcal{F}(2)_n||_{\infty}=  F_n$
    \item $||\mathcal{F}(2)_n||_{-\infty}=  0 $
    \item $||\mathcal{F}(2)_n||_{-1} = \infty$
    \item $||\mathcal{F}(2)_n||_0^n= 0$. \\
\end{enumerate}

A limitation of the symmetric $k$-Fibonacci matrix $\mathcal{Q}(k)_n=\left[q(k)_{ij} \right]$ described in \cite{lee2003} is that it can not be fully represented by $k$-Fibonacci numbers. Let us define a new $n \times n$ $\mathcal{S}$-type symmetric $k$-Fibonacci matrix, denoted by $\mathcal{S}(k)_n=\left[s(k)_{ij} \right]$, and defined as:
\begin{equation}
    s(k)_{ij}= s(k)_{ji} = g(k)_{i+j-1}.
\end{equation}

This definition ensures that all matrix elements are $k$-Fibonacci numbers and each element is the sum of $k$ previous elements in rows and columns. Then, a $2$-Fibonacci symmetric matrix of order $n$ is:
\begin{equation}
\mathcal{S}(2)_n = \left[ \begin{array}{ccccc} F_1 & F_2 & F_3 & \ldots & F_n \\ F_2 & F_3 & F_4 & \ldots & F_{n+1} \\ F_3 & F_4 & F_5 & \ldots & F_{n+2}\\   \vdots & \vdots & \vdots & \ddots & \vdots \\  F_n & F_{n+1} & F_{n+2} & \ldots & F_{2n-1} \end{array} \right].
\end{equation}

Reorganizing the elements of $\mathcal{S}(2)_n$ leads to a triangular structure with a clear pattern:
\begin{equation}
\mathcal{S}^*(2)_n = \left[ \begin{array}{ccccc} F_1 & 0 & 0 & \ldots & 0  \\ F_2 & F_2 & 0 & \ldots & 0 \\ F_3 & F_3 & F_3 & \ldots & 0\\   \vdots & \vdots & \vdots & \cdots & \vdots \\ F_{n-1} & F_{n-1} & \ldots & F_{n-1} & 0 \\  F_n & F_n &  \ldots & F_n & F_n \\ F_{n+1} & F_{n+1} & \ldots & F_{n+1} & 0 \\ \vdots & \vdots  & \vdots &  \cdots &  \vdots \\ F_{2n-2} & F_{2n-2} &  0& 0 & 0 \\ F_{2n-1} & 0 &  0& 0 & 0 \\ \end{array} \right].
\end{equation}

An example of this reorganization for the elements of $\mathcal{S}(2)_5$ is:
\begin{equation}
\begin{array}{ccccc} 1 &  &  &  &  \\ 1 & 1 &  &  &   \\ 2 & 2 & 2 &  & \\   3 & 3 & 3 & 3 &  \\  5 & 5 & 5 & 5 & 5 \\ 8 & 8 & 8 & 8 &  \\ 13 & 13 & 13 &  &  \\ 21 & 21 & & &   \\ 34 &  &  &  &   \end{array} 
\end{equation}

Note that the sum of the elements of the $i$-th row is $5-|5-i|$ times $F_i$. Generalizing to a matrix $\mathcal{S}(2)_n$ of order $n$, the sum of the elements of the $i$-th row is $n-|n-i|$ times $F_i$. This pattern leads to:
\begin{equation}
    ||\mathcal{S}(2)_n||_p^p = \sum_{i=1}^{2n-1} F_i^p(n-|n-i|).
\end{equation}

Setting $p$ to different values results in the following properties:

\begin{enumerate}
    \setlength\itemsep{1em}
    \setcounter{enumi}{27}
    \item $||\mathcal{S}(2)_n||_1 = \sum_{i=1}^{2n-1} F_i(n-|n-i|)$
    \item $||\mathcal{S}(2)_n||_2^2 = \sum_{i=1}^{2n-1} F_i^2(n-|n-i|)$
    \item $||\mathcal{S}(2)_n||_3^3 = \sum_{i=1}^{2n-1} F_i^3(n-|n-i|)$
    \item $||\mathcal{S}(2)_n||_{\infty}= F_{2n-1}$
    \item $||\mathcal{S}(2)_n||_{-\infty}= F_1 $
    \item $||\mathcal{S}(2)_n||_{-1} = \displaystyle{\left(\sum_{i=1}^{2n-1} \frac{(n-|n-i|)}{F_i}\right)^{-1}}$
    \item $||\mathcal{S}(2)_n||_0^n= \sum_{i=1}^{2n-1} F_i^{n-|n-i|}$. \\
\end{enumerate}


\section{Distances between vectors and matrices}

When considering two vectors (matrices) of the same size, the $p$-norm of the difference of the two vectors (matrices) implies a $p$-distance between them. To provide a general $p$-distance definition that applies to vectors and matrices, it is convenient to transform matrices into vectors through the following vectorization operator:

\begin{definition}
Given an $m\times n$ matrix $A$, vectorization operator $\mathrm{vec}(A)$ produces a vector $\boldsymbol{x}$ of $m \cdot n$ elements obtained by arranging the elements of $A$ in row-major order, i.e., by arranging them sequentially row by row:
\begin{equation}
    x_{i+j-1}=a_{ij} \; \forall i\in\{1,\ldots,m\} \; \forall j\in\{1,\ldots,n\}.
\end{equation}
\end{definition}

Then, the $p$-distance between two $n$-dimensional Fibonacci vectors $\boldsymbol{x}$ and $\boldsymbol{y}$ is defined as:

\begin{definition}
Given $\boldsymbol{x}, \boldsymbol{y} \in \mathbb{F}(k)^n$, where $\mathbb{F}(k)^n$ denotes the set of all $k$-Fibonacci vectors of size $n$, the $p$-distance between $\boldsymbol{x}$ and $\boldsymbol{y}$  is defined as
\begin{equation}
    ||\boldsymbol{x}-\boldsymbol{y}||_p = \left( \sum_{i=1}^n |x_i - y_i|^p\right)^{1/p}.
\end{equation}
\end{definition}

Consider two $(n,r)$-Fibonacci vectors whose elements are separated by $d$ positions in the sequence $\boldsymbol{x} = \boldsymbol{q}_{n+d,r}$ and $\boldsymbol{y} = \boldsymbol{q}_{n,r}$. Then:
\begin{equation}
    \boldsymbol{x} - \boldsymbol{y} = (F_{n+d+1}-F_{n+1}, F_{n+d+2}-F_{n+2}, \ldots, F_{n+d+r} - F_{n+r}) 
\end{equation}
and
\begin{equation}
    ||\boldsymbol{x}-\boldsymbol{y}||_p^p = (F_{n+d+1}-F_{n+1})^p + (F_{n+d+2}-F_{n+2})^p + \ldots +(F_{n+d+r} - F_{n+r})^p.
\end{equation}

From the fundamental identity of Fibonacci numbers $F_{n+1}-F_{n-1} = F_n$, for $d=2$ we have:
\begin{equation}
    ||\boldsymbol{x}-\boldsymbol{y}||_p^p (d=2) = F_{n+2}^p + F_{n+3}^p + \ldots +F_{n+r+1}^p = \sum_{i=1}^r F_{n+i+1}^p.
\end{equation}

Setting $d=3$ implies that $F_{n+4}-F_{n+1} = F_3 \cdot F_{n+2}$ and $F_{n+3+r}-F_{n+r} = F_3 \cdot F_{n+r+1}$, leading to: 
\begin{equation}
    ||\boldsymbol{x}-\boldsymbol{y}||_p^p (d=3) = (F_3 \cdot F_{n+2})^p + \ldots + (F_3 \cdot F_{n+r+1})^p = F_3^p \sum_{i=1}^r F_{n+i+1}^p.
\end{equation}

Consider now the sum and the difference of two $(n,r)$-Fibonacci vectors whose elements are separated by $d$ positions in the sequence:
\begin{equation}
    \boldsymbol{x} + \boldsymbol{y} = (F_{n+d+1}+F_{n+1}, F_{n+d+2}+F_{n+2}, \ldots, F_{n+d+r} + F_{n+r}) 
\end{equation}
\begin{equation}
    \boldsymbol{x} - \boldsymbol{y} = (F_{n+d+1}-F_{n+1}, F_{n+d+2}-F_{n+2}, \ldots, F_{n+d+r} - F_{n+r}).
\end{equation}

Then, the addition of the $1$-norms of the sum and difference vectors is:
\begin{equation}
    ||\boldsymbol{x}+\boldsymbol{y}|| + ||\boldsymbol{x}-\boldsymbol{y}|| = 2(F_{n+d+1}+ \ldots+F_{n+d+r}) = 2(F_{n+r+d+2} - F_{n+d+2})
\end{equation}

Similarly, using the parallelogram law in inner product spaces for $2$-norms, we have:
\begin{equation}
    ||\boldsymbol{x}+\boldsymbol{y}||_2^2 + ||\boldsymbol{x}-\boldsymbol{y}||_2^2 =  2(||\boldsymbol{x}||_2^2 + ||\boldsymbol{y}||_2^2) 
\end{equation}

\begin{equation}
    ||\boldsymbol{x}||_2^2 = F_{n+r+d} \cdot F_{n+r+d+1} - F_{n+d} \cdot F_{n+d+1}
\end{equation}

\begin{equation}
    ||\boldsymbol{y}||_2^2 = F_{n+r} \cdot F_{n+r+1} - F_{n} \cdot F_{n+1}
\end{equation}

\begin{equation}
    2(||\boldsymbol{x}||_2^2 + ||\boldsymbol{y}||_2^2) =  2(F_{n+r+d} \cdot F_{n+r+d+1} - F_{n+d} \cdot F_{n+d+1} +  F_{n+r} \cdot F_{n+r+1} - F_{n} \cdot F_{n+1}).
\end{equation}

By considering the golden ratio:
\begin{equation}
    \lim_{n \to \infty} \frac{F_{n+1}}{F_{n}} = \varphi
\end{equation}
we can approximate the product of two consecutive Fibonacci numbers through:
\begin{equation}
F_n \cdot F_{n+1} \approx \varphi  F_n^2    
\end{equation}
to further simplify the addition of the $2$-norms of the sum and difference of $(n,r)$-Fibonacci vectors whose elements are separated by $d$ positions in the sequence:
\begin{equation}
    ||\boldsymbol{x}+\boldsymbol{y}||_2^2 + ||\boldsymbol{x}-\boldsymbol{y}||_2^2=2(||\boldsymbol{x}||_2^2 + ||\boldsymbol{y}||_2^2) \approx 2 \varphi( F_{n+r+d}^2 -  F_{n+d}^2 +  F_{n+r}^2  - F_{n}^2).
\end{equation}

\section{Conclusions}

By defining a $p$-norm over the set of vectors and matrices populated with Fibonacci numbers, several essential properties are derived by varying parameter $p$. Given the decreasing character of the difference of $p$-norms as $p$ increases, a closed-form expression is given to obtain the value of $p$ that makes the difference between the $p$-norm and the infinite norm as small as needed. A new class of symmetric $k$-Fibonacci matrix is defined, and reorganizing its elements leads to a pattern that simplifies the computation of its $p$-norms. The analysis of $p$-distances between same-sized vectors (matrices) is illustrated by specific examples when the elements are separated by two and three positions in the sequence.

\medskip

\noindent MSC2020: 11B39, 15A60 

\end{document}